\input amstex
\documentstyle{ams-j}
\NoBlackBoxes

\define\inc{\subseteq}

\define\m{\bold m}

\define\la{\lambda}
\define\Iq{I^{[q]}}

\define\h{\operatorname{height}}
\define\di{\operatorname{dim}}

\topmatter
\title Unmixed Local Rings With Minimal Hilbert-Kunz Multiplicity Are Regular\endtitle
\rightheadtext{Hilbert-Kunz Multiplicity}
\author Craig Huneke and Yongwei Yao\endauthor

\address
Department of Mathematics, University of Kansas,
Lawrence, KS 66045
\endaddress
\email
huneke\@math.ukans.edu
\endemail
\address
Department of Mathematics, University of Kansas,
Lawrence, KS 66045
\endaddress
\email
yyao\@math.ukans.edu
\endemail

\issueinfo{00}
{0}
{Xxxx}
{2000}

\thanks The first author was partially supported by the NSF.\endthanks

\abstract
We give a new and simple proof that unmixed local rings having
Hilbert-Kunz multiplicity equal to $1$ must be regular.
\endabstract

\endtopmatter

\document

\head 1. Introduction \endhead

In 1969, Kunz \cite{Ku1} introduced a new numerical function of a Noetherian local
ring in positive characteristic.
One of Kunz's main results was the following characterization of regular
local rings in positive characteristic $p$. We always let $q$ denote a varying power
of $p$, and for an ideal $I$ write $\Iq = (i^q|\, i\in I)R$.

\proclaim{Theorem 1.1 \cite{Ku1, Prop. 3.2, Thm. 3.3}} Let $(R,\m)$ be a Noetherian local ring
of positive characteristic $p > 0$. Then 
\roster
\item $\la_R(R/\m^{[q]})\geq q^d$ for all $q = p^e$ ($e\geq 1$).
\item $R$ is regular iff for some $q = p^e$ with $e\geq 1$ (equivalently for all $q$),
$\la_R(R/\m^{[q]}) = q^d$.
\endroster
\endproclaim

In \cite{Ku2} Kunz observed that the function $q:\rightarrow  \la_R(R/\m^{[q]})/q^d$ should
give interesting information about the nature of the singularity of $R$. 
His interest was at least partly inspired by searching for numerical
invariants of singularities which might behave well under blowing up,
as a means to understand resolution of singularities in positive 
characteristic. The idea of taking the limit of these numbers as $q$
goes to infinity was discussed by Kunz, and that the limit exists was
shown by Monsky. He introduced the following definition:

\definition{Definition 1.2 \cite{Mo}} Let $(R,\m)$ be a Noetherian local ring
of dimension $d$ and positive prime characteristic $p$.
Let $M$ be a finitely generated $R$-module, and $I$ an $\m$-primary ideal.
Define $$ e_{HK}(I,M) = \underset q\to \varinjlim \frac {\la(M/\Iq M)} {q^d}.$$
\enddefinition

If the dimension of $M$ is strictly less than the dimension of $R$, this
limit is $0$. When dim$(M) = \di(R)$,  the limit always exists and is a positive real number
(Monsky \cite{Mo}). By definition the Hilbert-Kunz multiplicity of $R$ is
$e_{HK}(R) = e_{HK}(\m,R)$. We write $e_{HK}(I) = e_{HK}(I,R)$. Theorem
1.1 above proves that $e_{HK}(R)\geq 1$, and if $R$ is regular, $e_{HK}(R) = 1$.
 This led to a natural question,
probably first posed by Kunz in conversations: 

\remark{Question 1.3} If $(R,\m)$ is a Noetherian unmixed local
ring of positive characteristic, and $e_{HK}(R) = 1$, is $R$ regular?
\endremark 
\smallskip
Recall that a local ring $R$ is \it unmixed \rm if 
$\di(\hat R) = \di(\hat R/Q)$ for every associated prime $Q$ of
$\hat R$. The same assumption is
needed in the famous result of Samuel (in the case where $R$ contains
a field \cite{Sa}) and Nagata (in the general case) that the usual multiplicity of
an unmixed local Noetherian ring is $1$ iff the ring is regular.

\proclaim{Theorem 1.4 \cite{Na, Thm. 40.6}} Let $R$ be an unmixed Noetherian local ring, not necessarily
of characteristic $p>0.$ 
If $R$ has multiplicity one, then $R$ is regular.
\endproclaim
 
In recent beautiful and intricate work,  K.-I. Watanabe and K.
Yoshida \cite{WY, Theorem 1.5} have shown that Question 1.3 has a positive answer. 
In this paper we give a new proof, which 
avoids some of the harder parts in the proof of Watanabe and
Yoshida; in particular we avoid reference to standard systems of parameters
and even to the theory of tight closure, although the present proof was
inspired through considerations involving tight closure.

When $R$ is regular, the Hilbert-Kunz multiplicity is easy to compute,
due to the exactness of Frobenius. One can easily prove (see \cite{BH, Ex. 8.2.10}):

\proclaim{Proposition 1.5} Let $(R,\m)$ be a regular local Noetherian ring of prime
characteristic $p> 0$ and dimension $d$. Let $I$ is a $m$-primary ideal. Then
for all $q$,
$$\lambda(R/I^{[q]}) = q^d \cdot \lambda (R/I). \tag{1.5.1}$$ 
In particular, $e_{HK}(I) = \lambda (R/I)$ for all $\m$-primary ideals $I$.
\endproclaim

\head 2.  Preliminaries \endhead

The next lemma and its corollaries 
are very useful in studying the length of $R/I^{[q]}$
for a $\m$-primary ideal. They can be found in \cite{WY, Proposition 4.1, 
Lemma 4.2} in slightly different forms. The filtration argument used in
the proof can also be found in [Ha, Proposition 5.2.1]. We include the
proof here for completeness.

\proclaim{Lemma 2.1} Let $(R,\m)$ be a Noetherian local ring of
characteristic $p>0$. Let $I \subseteq J$ be two ideals with $I$ 
$\m$-primary (we allow $J = R$). Then
$$
\lambda (R/I^{[q]}) \leq  \lambda (J/I) \cdot \lambda (R/\m^{[q]})+
\lambda(R/J^{[q]}).
$$
\endproclaim

\demo{Proof} Set $s = \lambda (J/I)$. Take a filtration of $I
\subseteq J \subseteq R$ 
$$
I=J_0 \subsetneq J_1 \subsetneq J_2 \subsetneq \cdots  \subsetneq J_s
= J  \subseteq R
$$
so that $\lambda (J_i/J_{i-1}) =1$ i.e. $J_i/J_{i-1} \cong R/\m, \
\forall i=1,2,\dots,s.$ That is to say $J_i = (J_{i-1}, x_i)$ for some
$x_i \in J_i$ such that $J_{i-1}:x_i = \m.$

For every $q=p^e,$ there is a corresponding filtration of $I^{[q]}
\subseteq J^{[q]} \subseteq R$
$$
I^{[q]}=J_0^{[q]} \subseteq J_1^{[q]} \subseteq J_2^{[q]} \subseteq
\cdots  \subseteq J_s^{[q]} = J^{[q]}  \subseteq R,
$$
where $J_i^{[q]}/J_{i-1}^{[q]} \cong R/(J_{i-1}^{[q]}: x_i^q)$, which is
a homomorphic image of $R/m^{[q]},$ for every $i=1,2,\dots, s.$ So 
$\lambda (J_i^{[q]}/J_{i-1}^{[q]}) \leq \lambda (R/\m^{[q]}).$ Therefore
$\lambda (R/I^{[q]}) \leq  \lambda (J/I) \cdot \lambda (R/\m^{[q]})+
\lambda(R/J^{[q]}).$ \qed
\enddemo

\proclaim{Corollary 2.2}  Let $(R,\m)$ be a Noetherian local ring of
characteristic $p>0$. Let $I$ be a $\m$-primary ideal of $R$. Then
\roster
\item
$\lambda (R/I^{[q]}) \leq  \lambda (R/I) \cdot \lambda (R/\m^{[q]}).$
\item If $I\inc J$ then $e_{HK}(I,R)\leq \la(J/I)e_{HK}(R) + e_{HK}(J,R)$.
\endroster

\endproclaim

\demo{Proof} In case (1), we take $J = R$ and apply Lemma 2.1.
In case (2) the Corollary  follows from Lemma 2.1 by 
dividing by $q^{\dim R}$ and then taking the limits. \qed

\enddemo

An important ingredient in our proof is a calculation which shows that
for some ideals, the Hilbert-Kunz multiplicity is well-behaved.

\proclaim {Theorem 2.3} Let $(R,\m)$ be a Noetherian local ring of
characteristic $p>0$ and dimension $d$. Let $J$ be an ideal such that $\dim R/J =1$ and $\h J = d-1$.
Assume that $x\in R$ is a  non-zerodivisor
in $R/J$, and set $I = (J,x)$. Assume that $R_P$ is regular for every minimal
prime $P$ above $J$. Then
$$
e_{HK}\left(I, R\right) \geq \lambda (R/I).
$$
\endproclaim
 
\demo{Proof} Using the properties of the usual multiplicity of
parameter ideals, the associativity formula for the usual multiplicity,
and (1.5.1), we have
$$ \align
e_{HK}(I,R) &= \lim_{q \to \infty}\frac 1{q^d} \cdot \lambda
(R/I^{[q]}) =
\lim_{q \to \infty}\frac 1{q^d} \cdot \lambda (R/(J^{[q]},x^q)) \\
&\geq \lim_{q \to \infty}\frac 1{q^d} \cdot e(x^q; R/J^{[q]}) =
\lim_{q \to \infty}\frac q{q^d} \cdot e(x; R/J^{[q]}) =
\lim_{q \to \infty}\frac 1{q^{d-1}} \cdot e(x; R/J^{[q]})\\
&=\lim_{q \to \infty}\frac 1{q^{d-1}} \cdot \sum_{P \in \min(R/J)}
e(x;R/P) \cdot \lambda_{R_P}(R_P/J_P^{[q]})\\
&=\lim_{q \to \infty}\frac 1{q^{d-1}} \cdot \sum_{P \in \min(R/J)}
e(x;R/P) \cdot q^{d-1} \cdot \lambda_{R_P}(R_P/J_P) \\
&=\lim_{q \to \infty} \sum_{P \in \min(R/J)}
e(x;R/P) \cdot \lambda_{R_P}(R_P/J_P)\\
&=\sum_{P \in \min(R/J)}e(x;R/P) \cdot \lambda_{R_P}(R_P/J_P)=
e(x;R/J)=\lambda (R/(J,x))\\
&= \lambda (R/I). 
\endalign  $$\qed
\enddemo

\remark{Remark 2.4} 
One might suspect equality in every case in Theorem 2.3. However
the following example provided by the referee shows that inequality
is the best one can do:
Let $R=k[[x,y,z]]/(xy-z^n), J=(y,z), I=m$.  Then
$e_{HK}(I)=(2n-1)/n > \lambda(R/I)=1$.

However, if in addition we assume that $e_{HK}(R) = 1$, then
it follows from Theorem 2.3 and Corollary 2.2 that in fact
$$e_{HK}\left(I, R\right) = \lambda (R/I).$$
However, we will not use this equality in the sequel.
In any case, after we prove that $e_{HK}(R) = 1$ implies
the regularity of $R$, $e_{HK}(I) = \lambda(R/I)$ for
all $m$-primary ideals $I$.

\endremark

\bigskip
 \head3. A Criterion for Regular Rings  \endhead

The critical step in proving our main result is in constructing an
$\m$-primary ideal $I\inc \m^{[p]}$ such that $e_{HK}(I)\geq \la(R/I)$.
In the paper of Watanabe and Yoshida \cite{WY}, this step also
played an important role. Their construction was
done by taking $I$ to be an ideal generated by parameters. However, to
prove the inequality it was necessary for them to first prove the
ring $R$ is forced to be Cohen-Macaulay if its Hilbert-Kunz multiplicity is
one; this is the difficult part of
their proof, and required tools from the theory of tight closure and
work on standard systems of parameters. We are able to entirely avoid
this point by focusing our attention on ideals $I$ which are not
necessarily generated by parameters. Here is the theorem of Watanabe and Yoshida \cite{WY,
Thm. 1.5}.

\proclaim {Theorem 3.1} Let $(R,\m)$ be an unmixed
Noetherian
local ring of characteristic $p>0.$ 
If $e_{HK}(R)=1$, then $R$ is regular.
\endproclaim

\demo {Proof} 
Since the Hilbert-Kunz multiplicity
of $R$ is the same as that of its completion, we may assume $R$ is
complete.  The associativity formula for Hilbert-Kunz multiplicity shows
that
$$e_{HK}(R) = \sum _P e_{HK}(R/P)\cdot \la(R_P)$$
where the sum is over all minimal primes of maximal dimension.
Since $e_{HK}(R) = 1$, we deduce that $R$ can have only one minimal
prime $P$ and $R_P$ has to be field, i.e. $P_P=0.$ Hence $P=0$ since
$R\setminus P$ consists of non-zero divisors. Thus $R$ is a domain.

It suffices to prove that $\la(R/\m^{[p]})\leq p^d$ (where $d = \di(R)$) as then Theorem 1.1 (1)
first gives equality, and then Theorem 1.1 (2) gives that $R$ must be
regular.

The singular locus of $R$ is closed and not equal to Spec$(R)$. It follows
we can choose a prime $P$ such that dim$(R/P) = 1$ and $R_P$ is regular.  
(If not, the intersection of all primes $Q$ such that dim$(R/Q) = 1$ must
be nonzero. Choose $f\ne 0$ in this intersection and extend $f$ to a
full system of parameters $g_1,...,g_{d-1},f$. Taking a minimal prime $Q$
over $(g_1,...,g_{d-1})$ gives a contradiction.) 
Since the intersection of the symbolic powers of $P$ is zero and $R$ is
complete, Chevalley's lemma gives that some sufficiently large
symbolic power of $P$ lies inside $\m^{[p]}$. Call this symbolic power $J$.
Choose $x\in \m^{[p]}$ such that $x\notin P$. The ideal $I = (J,x)$ lies
in $\m^{[p]}$ and satisfies the hypothesis of Theorem  2.3. Hence
$$e_{HK}(I) \geq \la(R/I).$$

On the other hand we have $e_{HK}(I,R) \leq \lambda(\m^{[p]}/I) \cdot
e_{HK}(R) + e_{HK}(\m^{[p]},R) = \lambda(\m^{[p]}/I) +
e_{HK}(\m^{[p]},R) \leq \lambda(\m^{[p]}/I) + \lambda(R/\m^{[p]}),$ by
Lemma 2.1 and Corollary 2.2.
 
That is to say  
$$ \align 
\lambda(\m^{[p]}/I) + \lambda(R/\m^{[p]}) =&\lambda (R/I) \leq e_{HK}(I,R) \\
\leq &\lambda(\m^{[p]}/I) + e_{HK}(\m^{[p]},R) \\
\leq &\lambda(\m^{[p]}/I) + \lambda(R/\m^{[p]}),
\endalign $$
which forces $\lambda(R/\m^{[p]}) = e_{HK}(\m^{[p]},R)$. However,
$$e_{HK}(\m^{[p]},R) = \varinjlim \frac {\la(R/\m^{[pq]})} {q^d} =
\varinjlim \frac {p^d\cdot\la(R/\m^{[pq]})} {(pq)^d} = p^d\cdot e_{HK}(R) = p^d.$$ 
Together the equalities  imply that $\lambda(R/\m^{[p]}) = p^d$, which
implies that $R$ is regular by Theorem 1.1. \qed   

\enddemo

\remark{Remark 3.2}
An alternate proof could be given which by induction allows one to assume that
$R_Q$ is regular for all primes $Q\ne \m$. This is due to the upper semi-continuity of the
Hilbert-Kunz multiplicity. This was shown by Kunz \cite{Ku2}.
We append a shorter proof, much in the spirit of our simplification 
of the main theorem.
\endremark
 
\proclaim{Theorem 3.3 \cite{Ku2, Cor. 3.8}}
Let $(R, \m)$ be a Noetherian local ring of
characteristic $p>0$, and let P be a prime ideal of R such that
$\h(P) + \dim(R/P) = \dim(R)$. Then $e_{HK}(R_P)\leq e_{HK}(R).$  In
fact,  if $t = \dim(R/P)$, then  $q^t \cdot ~\lambda_{R_P}((R/P^{[q]})_P)
\leq \lambda(R/\m^{[q]})$ for  every  $q=p^e.$
\endproclaim
 
\demo{Proof} By induction, it is enough to prove the case where $
ht(P) =\di(R) - 1$. Notice it suffices to prove the second inequality.
 
Choose $f \in \m - P$. Then, using the properties of the usual
multiplicity of parameter ideals, the associativity formula for the
usual multiplicity, we have,  for all $q=p^e$,
$$\align
\lambda(R/(P, f)^{[q]}) &= \lambda(R/(P^{[q]}, f^q))\\
& \geq e(f^q; R/P^{[q]})\\
& =\lambda _{R_P}((R/P^{[q]})_P) \cdot e(f^q; R/P)\\
& =\lambda _{R_P}((R/P^{[q]})_P) \cdot q  \cdot \lambda(R/(f, P)).
\endalign$$
Also, by Corollary 2.2, we know  that
$\lambda(R/(f, P)) \cdot \lambda(R/\m^{[q]})\geq \lambda(R/(P, f)^{[q]}).$
Hence  $\lambda(R/\m^{[q]}) \geq q \cdot \lambda_{R_P}((R/P^{[q]})_P)$
for every $q=p^e.$  \qed
\enddemo

\bigskip

\centerline{\bf Bibliography}
\bigskip
\refstyle{A}

\enddocument